\documentclass{article}

\begin{document}

\title{Axiomatic Origins of Mathematical Entropy: Grading Ordered Sets}

\author{Alexander Dukhovny\\
\small{Department of Mathematics, San Francisco State University} \\
\small{San Francisco, CA 94132, USA} \\
\small{dukhovny[at]sfsu.edu}\\
}

\date{\today}          

\maketitle

\begin{abstract}

Shannon's entropy and other entropy-based concepts are derived from the new, more general concept of relative divergence of one ``grading'' function on a linearly  ordered set from another such function. The definition of relative divergence is derived based on ``common sense'' assumptions about comparing grading functions. Shannon's entropy formulas emerge from the respective relative divergence ones, entropy based methods are  extended to more general cases and some new applications.

\end{abstract}

\section{Introduction}
\label{sec:Intro}

Shannon's entropy is one of the most important and effective concepts in mathematics. It has been used in many ways, especially in the framework of the Maximum Entropy Principle (MEP) serving as a mathematical representation of the Insufficient Reason Principle, to select the ``least assuming'' probability distribution under constraints. The references are too many to quote.

There have been numerous generalizations of the original Shannon's entropy definition: relative entropy, Kullback-Leibler divergence, partition entropy, Kolmogorov-Sinai entropy, topological entropy, entropy of general non-probabilistic measures (capacities) and great many others.

The concept itself was introduced by Shannon to quantify the information contained in a coded signal \cite{Shannon}. The effectiveness of Shannon's entropy led to its numerous  applications in other fields of mathematics. In the process, the axioms leading to the Shannon's entropy formula were established (see, e.g., \cite{Jaynes}) and carried over from its original context of the probability theory to new applications. However, in all such works the starting point was the original Shannon's entropy formula, and no attempt has been made, as far as we searched, to derive it from a more general, more intuitive concept. 

The goal of this paper is precisely that: to introduce a natural mathematical concept arising within a very basic context: assigning real-valued grades to elements of a linearly  ordered set, increasing in agreement with the order. Using some simple principles, we introduce the concept of ``relative divergence of one grading function from another'' (using the Kullback-Leibler terminology (see \cite{K-L})).  In turn, relative divergence formulas yield Shannon's entropy and other related formulas for well-known cases and open new ways to look at new applications.

The general idea of the proposed approach begins as follows: let $W$ be a linearly ordered set and use $\prec$ to denote the ordering relation of its elements. A real-valued function $F : W \to R$ is said to be a grading function on $W$ if 

$w \prec v \iff F(w) < F(v)$ for all $w, v \in W $.

(For example, if $W$ is enumerated then the position-function of an element is a natural grading function on $W$.) 

As such, a grading function's values can be used to define the ordering relation $\prec$ on $W$. Also, its inverse function $w = F^{-1}(u)$is defined  for all $u \in im(F) \subset R$.

Suppose two grading functions are defined on $W \colon$
$u = F(w)\colon W \rightarrow R$ and $G(w)\colon W \rightarrow R$. 

The definition of relative divergence of $F$ from $G$ over the set $W$ is derived here in Section 4 from some natural assumptions and is presented by two conceptually equivalent formulas depending on the nature of $W$ and $im(F)$.

When $W = \{ \ldots , w_{-1}, w_0 , w_1 , \ldots \}$ is a countable ordered set,  relative divergence of $F$ from $G$ on $W$ is defined as

\begin{equation} \label{eq:(1)}
\mathcal{D}(F \Vert G) = \ \sum_{k= - \infty}^\infty\
\ln  \left( \frac{\Delta_k G}{\Delta_k F }\right) \Delta_k F , 
\end{equation}

$\Delta_k F = F(w_k) - F(w_{k-1}), \hspace{2mm} \Delta_k G= G(w_k) - G(w_{k-1}), \hspace{2mm} k = \ldots,-1,0,1,\ldots .$

When $im(F)$ is an interval, relative divergence of $F$ from $G$ on $W$ is defined as

\begin{equation}\label{eq:(2)}
\mathcal{D}(F \Vert G)=
\int_{im(F)} \ln \left( \frac{d}{du}(G(F^{-1}(u))) \right)du
\end{equation}

provided that the integral converges absolutely.

In section 2 we use equation (1) to derive relative divergence formulas for countable linearly ordered sets arising in the context of the probability theory and measure theory. We obtain formulas for some important examples: classic Shannon's entropy, relative entropy, partition entropy and entropy of non-additive measures.

Similarly, in section 3 continuous cases are considered. The relative divergence formula (2) is shown to yield classical relative entropy formula for continuous probability distributions, and to correct the ``individual'' entropy formula in order to make the ``relative''and ``individual'' ones agree.

\section{Countable linearly ordered sets}
\label{sec:Countable}

In this section we consider cases where the linearly ordered set $W$ consists of nested subsets of another countable set $E$, and the ordering relation $\prec$ is the subset inclusion relation $\subset$.

\subsection{Probability theory applications}
\label{subsection: ProbsDiscrete}

In the context of the probability theory, suppose 
$E = \{ e_1 , e_2 , \ldots \}$ is a countable set. Enumerating its elements, let
$W$ be a sequence of nested subsets of $E$ as follows:
$w_0 = \emptyset, w_k = w_{k-1} \cup \{e_k \}, k=1, 2, \ldots.$

Suppose probability distributions $\{f_i\}$ and $\{g_i\}$ are defined on $E$: 
 
$$\sum_{i=1}^\infty f_i = 1, \hspace{3mm} \sum_{i=1}^\infty g_i = 1, \hspace{3mm} f_i , \hspace{2mm} g_i \geq 0, \hspace{2mm} \forall i.$$

Let 
$F(w_k) = \sum_{i=1}^k f_i, \hspace{2mm} G(w_k) = \sum_{i=1}^k g_i, \hspace{2mm} k=0,1,\ldots,$

then
$\Delta_k F = f_k, \hspace{2mm} \Delta_k G = g_k, \hspace{2mm} k=1,2,\ldots$

and equation (\ref{eq:(1)}) yields

\begin{equation}\label{eq:(3)}
\mathcal{D}(F \Vert G) = \ \sum_{k=1}^\infty\ f_k\ln \left(\frac{g_k}{f_k} \right), 
\end{equation}

which is the classical formula of relative entropy of the probability distribution $\{f_i\}$ w.r.t. the probability distribution $\{g_i\}.$

On the other hand, if $G(w_k) = k, \hspace{3mm}k = 0,1,\ldots,$

then 
$\Delta_k F= f_k, \hspace{3mm} \Delta_k G = 1, \hspace{3mm} k=1,2,\ldots,$ 

and equation (\ref{eq:(1)}) yields

\begin{equation}\label{eq:(4)}
\mathcal{D}(F \Vert G)  = \ \sum_{k=1}^\infty\  f_k \ln \left(\frac{1}{f_k } \right)
= - \ \sum_{k=1}^\infty\  f_k\ln({f_k }), 
\end{equation}
the Shannon's entropy of the probability distribution $\{f_i\}$.

\subsection{Measure theory applications}
\label{subsection: M-theory}

In the context of general measure theory, Shannon's entropy has been generalized to study new applications of interest (such as studies of ``capacities''), see, e.g., \cite{Dukh}, \cite{KMR}, \cite{Marichal}, \cite{KM}, \cite{HG}. The concept of relative divergence can be applied to study more general cases of such applications.

Let $\mu$ be a general measure on a space $S \colon$ a set-monotonic nonnegative function on $2^S$, not presumed additive or normalized (that is, that $\mu (S) = 1$). Let $\mathcal{P} = \{A_1, A_2, \ldots \}$ be a countable set of $\mu$-measurable non-overlapping subsets of $S$ making its partition. Enumerating the partition subsets, define $W$ as a sequence of nested subsets:
$w_0 = \emptyset,\hspace{3mm} w_k = w_{k-1} \cup \{A_k \}, \hspace{3mm} k=1, 2, \ldots.$

Let 
$F(w_k) = \sum_{i=1}^k \mu(A_i), \hspace{3mm} G(w_k) = k, \hspace{3mm} k=0,1,\ldots .$

Then
$\Delta_k F = \mu(A_k), \hspace{3mm} \Delta_k G = 1, \hspace{3mm} k=1,2,\ldots,$
and Equation (\ref{eq:(1)}) yields

\begin{equation}\label{eq:(5)}
\mathcal{D}(F \Vert G)  = - \ \sum_{k=1}^\infty\  \mu(A_k)\ln({\mu(A_k })), 
\end{equation}

which is known as the partition entropy (see, e.g., \cite{Sinai} for the case where  measure $\mu$ is additive and normalized).

Representing Shannon's entropy as relative divergence of a grading function on a countable linearly ordered set from the element's position function allows to extend the concept of entropy to the following case (see \cite{Dukh}).

Let $\mu$ be a general measure on a countable set $S$, and let $C = \{C_0 = \emptyset, \hspace{3mm} C_1, C_2, \ldots $ be a maximal chain of $\mu$- measurable subsets of $S.$ $C$ being an ordered set, $\mu_C(C_k)= \mu (C_k), \hspace{2mm} k=0, 1, \ldots$ is a grading function on $C$. Defining $G_C(C_k)= k, \hspace{2mm} k=0,1,\ldots$,  relative divergence of $\mu_C$ from $G_C$ is given (assuming absolute convergence of the series) by

\begin{equation}\label{eq:(6)}
\mathcal{D}(\mu_C \Vert G_C)
  = - \ \sum_{k=1}^\infty\  \Delta_k (\mu_C)\ln({\Delta_k (\mu_C). }) 
\end{equation}

Following \cite{Dukh}, in order to make sure that using the Maximum Entropy Principle will maximize all $\mathcal{D}(\mu_C \Vert G_C)$, the entropy of the general measure $\mu$ should be now defined as

\begin{equation}\label {eq:(7)}
\mathcal{H}(\mu) = \inf_C \mathcal{D}(\mu_C \Vert G_C).
\end{equation}

\section{Continuous linearly ordered sets}
\label{section: Continuous}

In cases where $W$ is not countable, the formula of equation (\ref{eq:(2)}) will apply as stated under the assumptions that $im(F)$ is an interval where the function $G(F^{-1}(u)$ is continuously differentiable and the integral in the equation converges absolutely.

In the context of the probability theory, we will consider continuous probability distributions and treat their cumulative distribution functions (c.d.f.) as grading functions on the set of intervals with increasing upper endpoints.

\subsection{Finite random variables}
\label{subsection: intervals}

Consider random variables $Y, X$ with the same set of possible values $E = [a,b]$, whose c.d.f.s $F(x),\hspace{2mm} G(x), \hspace{2mm}x \in E$ are continuously differentiable on $E$ with probability density functions $f(x), \hspace{2mm} g(x)$, respectively.

Let $W = \{ [a,x], \hspace{2mm} x \in E \}$. Treating $F(x), \hspace{2mm} G(x)$ as grading functions on $W, \hspace{2mm} im(F) = [0,1], \hspace{2mm} du = f(x)dx$,

$$\frac{d}{du}(G(F^{-1}(u))) = \frac{g(x)}{f(x}, \hspace{3mm} x=(F^{-1}(u)), \hspace{3mm} u\in E,$$

and equation (\ref{eq:(2)}) yields

\begin{equation}\label{eq:(8)}
\mathcal{D}(F \Vert G)=
\int_a^b f(x)\ln \left( \frac{g(x)}{f(x)} \right)dx,
\end{equation}

the classical formula of relative entropy of $Y$ w.r.t. $X$ in the continuous case.

When $X$ is uniformly distributed on $[a,b]$, then 
$$g(x) = \frac{1}{b-a}, \hspace{3mm} a\le x \le b,$$ 
so (\ref{eq:(8)}) yields

\begin{equation}\label{eq:(9)}
\mathcal{D}(F \Vert G)=
-\ln(b-a) - \int_a^b f(x)ln(f(x))dx = \mathcal{H}(Y) - ln(b-a),
\end{equation}

where $\mathcal{H}(Y))$ is the classical formula for the ``individual'' entropy of $Y$.

Now, following the definition of symmetric Kullback-Leibler divergence of probability distributions (see \cite{K-L}), symmetric divergence of two grading functions can be defined as

\begin{equation}\label{eq:(10)}
\mathcal{D}(F,G)= \mathcal{D}(F \Vert G)+\mathcal{D}(G \Vert F)= 
\end{equation}

$\int_{im(F)} ln(\frac{d}{du}(G(F^{-1}(u))))du +
\int_{im(G)} ln(\frac{d}{dv}(F(G^{-1}(v))))dv$

\section{Definitions and Proofs}
\label{DP}

Suppose a grading function $G$ is defined on a linearly ordered set $W$. If another grading function $F$ is introduced on $W$, we seek a general method to quantify the  information about the overall distinction of $F$ from $G$ on $W$, referred to as  ``relative divergence of $F$ from $G$ over $W$''. The intuitive way for that is to rate the change of grade made by $F \colon \hspace{2mm}  \Delta F = F(v) - F(w)$ 
in terms of the change of grade made by $G \colon \hspace{2mm}\Delta G = G(v) - G(w)$, 
where $w \prec v, \hspace{2mm} w,v \in W$. As a tool to quantify that information we seek a suitable function $h(s,t)$ to apply to $s = \Delta G, \hspace{2mm} t=\Delta F$. 

For that purpose $h(s,t)$ must be a continuously differentiable function, so that it would vary smoothly in response to infinitesimal deviations of its arguments. Also, to apply $h(s,t)$ to identical grading functions, we set $h(s,s)= 0$.

\textbf{Assumption 1}. When the same linear transformation is applied to both grading functions, the information on the relative divergence rate does not change, so

$$h(\Delta G, \Delta F)= h(\Delta (aG+b), \Delta (aF+b)), \hspace{3mm} \forall a, b \in R.$$

It follows trivially from Assumption 1 that $h$ has to be a function of the ratio of its arguments:

$$h(s,t) = h(s/t) = h(r_{GF})$$

where 
$$r_{GF} = \frac{\Delta G}{\Delta F}$$

is the rate of change of $G$ per unit change in $F$. 

In particular, $h(s,s)=h(1)=0.$

\textbf{Assumption 2}. For any grading functions $G, F, K$

$$h(r_{GK}) = h(r_{GF}) + h(r_{FK}).$$

Since, obviously, $r_{GK} = r_{GF}r_{FK}$, it follows directly that the only function $h(r), \hspace{2mm} x > 0,$ such that 
$$h(xy) = h(x)+h(y), \hspace{3mm} \forall x, y > 0, \hspace{3mm} h(1) = 0,$$
is (up to a constant factor) $h(x) = ln{x}$.

Choosing the constant factor to be 1, we obtain 
$$h(\Delta G, \Delta F)= \ln \left( \frac{\Delta G}{\Delta F} \right)$$

Now, $h(r_{GF})$ is relative divergence of $F$ from $G$ per unit change in $F$. To account for the relative divergence over the entire pair $\{w,v\}$ it is only natural to multiply that rate by the entire grade change of $F$ over that pair:

\textbf{Assumption 3}. Relative divergence of $F$ from $G$ over $W = \{w, v\}$

$\mathcal{D}(F \Vert G) = (\Delta F) h(\Delta G, \Delta F)$

\textbf{Assumption 4}. If $W$ is partitioned into two parts concatenated at a common element $v \in W$, then the entire relative divergence over $W$ is a sum of relative divergences over each part:

$$W = (W_1 \colon \{ \forall w \preceq v\} \cup W_2 \colon \{\forall w \succeq v \})
\Longrightarrow$$

$$\mathcal{D}(F\Vert G)|_W = \mathcal{D}(F\Vert G)|_{W_1} + \mathcal{D}(F \Vert G)|_{W_2}$$.

Now, the additivity property of Assumption 4 obviously extends to any concatenated partition, finite or countable (provided the emerging series converges absolutely), which gives proof to the formula of equation (\ref{eq:(1)}).

Denoting $u = F(w), \hspace{3mm} q(u) = G(F^{-1}(u))$, equation (\ref{eq:(2)}) can be rewritten as

\begin{equation} \label{eq:(11)}
\mathcal{D}(F \Vert G) = \ \sum_{k= - \infty}^\infty\
\ln \left( \frac{\Delta_k q}{\Delta_k u } \right) \Delta_k u, 
\end{equation}

where 
$$\Delta_k q = q(u_k) - q(u_{k-1}), \hspace{2mm} \Delta_k u = (u_k - u_{k-1}), \hspace{2mm} u_k = F(w_k), \hspace{2mm} k = \ldots,-1,0,1,\ldots.$$

Equation (\ref{eq:(11)}) now becomes a basis for the continuous case formula.

Let $W$ be a linearly ordered set such that $im(F)= [a,b]$, and assume that function 
$q(u)= G(F^{-1}(u))$ is continuously differentiable on $[a,b]$. 
Partition $[a,b]$ into 
$\mathcal{P} \colon \{u_0 = a, \hspace{2mm} u_1, \ldots, u_n = b \}$. 
Since $F$ (as a grading function) is invertible, consider a set $W(n) = \{w_k = F^{-1}(u_k), \hspace{2mm} k = 0,1, \ldots,n \}.$ By (\ref{eq:(11)}),

\begin{equation} \label{eq:(12)}
\mathcal{D}(F \Vert G)|_{W(n)} = \ \sum_{k=0}^n
\ln \left( \frac{\Delta_k q}{\Delta_k u } \right) \Delta_k u, 
\end{equation}

which is the Riemann sum of the function $\hspace{2mm} ln(\frac{d}{du}q(u))$. Since, by construction, $q(u)$ is positive-valued and is assumed to be continuously differentiable on $[a,b]$, as $\Vert \mathcal {P} \Vert \to 0$, the RHS of (12) approaches 

\begin{equation}\label{eq:(13)}
\mathcal{D}(F \Vert G)=
\int_a^b \ln ( \frac{d}{du}(q(u)))du
\end{equation}

which is the integral in equation (\ref{eq:(2)}) for the case where $im(F)=[a,b].$ 

Extending equation (\ref{eq:(13)}) to cases where $im(F)$ is an interval that is not closed at the endpoints, finite or infinite, in addition to assuming that $G(F^{-1})$ is continuously differentiable on $im(F)$, another condition would have to be satisfied by that function to guarantee absolute convergence of the emerging improper integral at both endpoints of $im(F)$.

\section{Conclusion} \label{end}

In this paper we introduced the concept of relative divergence of one grading function $F$ on a linearly ordered set $W$ from another grading function $G$. The formula for relative divergence in cases when $W$ is countable was presented in equation (\ref{eq:(1)}). Equation (\ref{eq:(2)}) expresses relative divergence for ``continuous'' cases where $im(F)$ is an interval. In section 4, the definition of relative divergence and both expressions are derived based on some natural assumptions about comparing grading functions.

In sections 2 and 3, in the context of the probability theory, the formulas for  relative divergence yield familiar formulas for Shannon's entropy and relative entropy for probability distributions on both discrete and continuous sample spaces. 

It turns out that Shannon's entropy in the discrete case emerges from relative divergence of the c.d.f. of the given distribution from the element's position function.

In the continuous case, to be consistent with relative entropy formula, the classical individual entropy formula must be corrected. Say, when the sample space of a random variable $Y$ is an interval $[a,b]$, equation (\ref{eq:(9)}) suggests that the individual entropy formula should be redefined as relative divergence of the c.d.f. of $Y$ from the c.d.f. of the uniform distribution on $[a,b] \colon $

\begin{equation}\label{eq:(13}
\mathcal{H}(Y) = \hspace{2mm} -\int_a^b f(x)ln((b-a)f(x))dx.
\end{equation}

Outside of the probability theory context, in section 2 the concept of relative divergence is applied to extend the known formulas for partition entropy and general measure entropy to cases where the measure is not presumed additive or normalized.

\section{Acknowledgements}

The author is deeply grateful to Sergei Ovchinnikov for his insightful comments and invaluable help in structuring the content and the flow of the paper.

\end{document}